\documentclass[11pt]{article}
\usepackage{amssymb}
\newtheorem{precor}{{\bf Corollary}}

\newenvironment{cor}{\begin{precor}{\hspace{-0.5
               em}{\bf.\ }}}{\end{precor}}
\newtheorem{precon}{{\bf Conjecture}}

\newenvironment{con}{\begin{precon}{\hspace{-0.5
               em}{\bf.\ }}}{\end{precon}}
\newtheorem{prealphcon}{{\bf Conjecture}}

\newenvironment{alphcon}{\begin{prealphcon}{\hspace{-0.5
               em}{\bf.\ }}}{\end{prealphcon}}
\newtheorem{predefin}{{\bf Definition}}

\newtheorem{preexm}{{\bf Example}}

\newtheorem{preappl}{{\bf Application}}

\newtheorem{prelem}{{\bf Lemma}}

\newtheorem{preproof}{{\bf Proof.\ }}

\newenvironment{proof}[1]{\begin{preproof}{\rm
               #1}\hfill{$\blacksquare$}}{\end{preproof}}
\newtheorem{prethm}{{\bf Theorem}}

\newenvironment{thm}{\begin{prethm}{\hspace{-0.5
               em}{\bf.\ }}}{\end{prethm}}
\newtheorem{prealphthm}{{\bf Theorem}}

\newenvironment{alphthm}{\begin{prealphthm}{\hspace{-0.5
               em}{\bf.\ }}}{\end{prealphthm}}

\newtheorem{prealphlem}{{\bf Lemma}}

\newenvironment{alphlem}{\begin{prealphlem}{\hspace{-0.5
               em}{\bf.\ }}}{\end{prealphlem}}

\newtheorem{prepro}{{\bf Proposition}}

\newtheorem{preprb}{{\bf Problem}}

\newtheorem{prerem}{{\bf Remark}}

\newenvironment{rem}{\begin{prerem}{\hspace{-0.5
               em}{\bf.\ }}}{\end{prerem}}
\newtheorem{preapp}{{\bf Application}}

\newtheorem{prequ}{{\bf Question}}

\def\conct[#1,#2]{\mbox {${#1} \leftrightarrow {#2}$}}
\def\dconct[#1,#2]{\mbox {${#1} \rightarrow {#2}$}}
\def\deg[#1,#2]{\mbox {$d_{_{#1}}(#2)$}}
\def\mindeg[#1]{\mbox {$\delta_{_{#1}}$}}
\def\maxdeg[#1]{\mbox {$\Delta_{_{#1}}$}}
\def\outdeg[#1,#2]{\mbox {$d_{_{#1}}^{^+}(#2)$}}
\def\minoutdeg[#1]{\mbox {$\delta_{_{#1}}^{^+}$}}
\def\maxoutdeg[#1]{\mbox {$\Delta_{_{#1}}^{^+}$}}
\def\indeg[#1,#2]{\mbox {$d_{_{#1}}^{^-}(#2)$}}
\def\minindeg[#1]{\mbox {$\delta_{_{#1}}^{^-}$}}
\def\maxindeg[#1]{\mbox {$\Delta_{_{#1}}^{^-}$}}
\def\isdef{\mbox {$\ \stackrel{\rm def}{=} \ $}}
\def\dre[#1,#2,#3]{\mbox {${\cal E}^{^{#3}}(#1,#2)$}}
\def\var[#1,#2]{\mbox {${\rm Var}_{_{#1}}(#2)$}}
\def\ls[#1]{\mbox {$\xi^{^{#1}}$}}
\def\hom[#1,#2]{\mbox {${\rm Hom}({#1},{#2})$}}
\def\onvhom[#1,#2]{\mbox {${\rm Hom^{v}}(#1,#2)$}}
\def\onehom[#1,#2]{\mbox {${\rm Hom^{e}}(#1,#2)$}}
\def\core[#1]{\mbox {$#1^{^{\bullet}}$}}
\def\cay[#1,#2]{\mbox {${\rm Cay}({#1},{#2})$}}
\def\sch[#1,#2,#3]{\mbox {${\rm Sch}({#1},{#2},{#3})$}}
\def\cays[#1,#2]{\mbox {${\rm Cay_{s}}({#1},{#2})$}}
\def\dirc[#1]{\mbox {$\stackrel{\rightarrow}{C}_{_{#1}}$}}
\def\cycl[#1]{\mbox {${\bf Z}_{_{#1}}$}}

\begin{document}
\begin{center}
{\Large \bf A Generalization of Kneser's Conjecture}\\
\vspace{0.3 cm}
{\bf Hossein Hajiabolhassan}\\
{\it Department of Mathematical Sciences}\\
{\it Shahid Beheshti University, G.C.}\\
{\it P.O. Box {\rm 1983963113}, Tehran, Iran}\\
{\tt hhaji@sbu.ac.ir}\\ \ \\
\end{center}
\begin{abstract}
\noindent We investigate some coloring properties of Kneser
graphs. A star-free coloring is a  proper coloring
$c:V(G)\rightarrow \Bbb{N}$ such that no path with three vertices
may be colored with just two consecutive numbers. The minimum
positive integer $t$ for which there exists a star-free coloring
$c: V(G) \rightarrow \{1,2,\ldots ,t\}$ is called the star-free
chromatic number of $G$ and denoted by $\chi_s(G)$. In view of
Tucker-Ky Fan's lemma, we show that $\chi_s({\rm
KG}(n,k))=2\chi({\rm KG}(n,k))-2=2n-4k+2$ provided that $n \leq
{8\over 3}k$. This gives a partial answer to a conjecture of
\cite{OMPO1}. Moreover, we show that for any Kneser graph ${\rm
KG}(n,k)$ we have $\chi_s({\rm KG}(n,k))\geq \max\{2\chi({\rm
KG}(n,k))-10, \chi({\rm KG}(n,k))\}$ where $n\geq 2k \geq 4$.
Also, we conjecture that for any positive integers $n\geq 2k \geq
4$ we have $\chi_s({\rm KG}(n,k))= 2\chi({\rm KG}(n,k))-2$.\\

\end{abstract}
\section{Introduction}
The local coloring of a graph $G$, defined in \cite{CHFA} and
\cite{CHSA}, is a proper coloring $c:V(G)\rightarrow \Bbb{N}$
such that no path with three vertices and no triangle may be
colored with just two and three consecutive numbers respectively.
In other words, for any set $S\subseteq V(G)$ with
$2\leq|S|\leq3$, there exist two vertices $u,v\in S$ such that
$|c(u)-c(v)|\geq m_S$, where $m_S$ is the number of edges of the
induced subgraph $G[S]$. The maximum color assigned by a local
coloring $c$ to a vertex of $G$ is called the value of $c$ and is
denoted by $\chi_l(G, c)$.  The local chromatic number of $G$ is
$\chi_l(G)=\displaystyle \min_c \chi_l(G,c)$, where the minimum
is taken over all local colorings $c$ of $G$. If
$\chi_l(G)=\chi_l(G,c)$, then $c$ is called a minimum local
coloring of $G$.

A star-free coloring is a  proper coloring $c:V(G)\rightarrow
\Bbb{N}$ such that no path with three vertices may be colored
with just two consecutive numbers. The minimum natural number $t$
for which there exists a star-free coloring $c: V(G) \rightarrow
\{1,2,\ldots ,t\}$ is called the star-free chromatic number of
$G$ and denoted by $\chi_s(G)$. Also, if $c: V(G) \rightarrow
\{1,2,\ldots ,t\}$ is a star-free coloring of $G$ and
$t=\chi_s(G)$, then $c$ is called a minimum star-free coloring.
Note that a local coloring is a star-free coloring with one more
requirement (no triangle may be colored with just three
consecutive numbers).

Hereafter, the symbol $[n]$ stands for the set $\{1,\ldots, n\}$.
Assume that $n \geq 2k$. The {\it Kneser graph} ${\rm KG}(n,k)$
is the graph with vertex set ${[n] \choose k}$, in which $A$ is
connected to $B$ if and only if $A \cap B = \emptyset$. For a
subset $X\subseteq [n]$ denote by ${X \choose k}$ the subgraph
induced by the collection of all $k$-subsets of $X$ in ${\rm
KG}(n,k)$. It was conjectured by Kneser \cite{KNE} in 1955, and
proved by Lov\'{a}sz \cite{LOV} in 1978, that $\chi({\rm
KG}(n,k))=n-2k+2$.

Let $c:V(G)\rightarrow \Bbb{N}$ be a proper coloring of a graph
$G$. If we replace the color $i$ by $2i-1$ for any positive
integer $1\leq i \leq \chi(G)$, then we obtain a local coloring
for G.

\begin{alphlem}{\rm \cite{CHSA}}\label{BOUND}
For any graph $G$ we have
$$\chi(G) \leq \chi_l(G)\leq 2\chi(G)-1.$$
\end{alphlem}

In view of definition of star-free chromatic number, one can
deduce that $\chi(G) \leq \chi_s(G) \leq \chi_l(G)\leq
2\chi(G)-1.$

In \cite{OMPO1}, the local chromatic number of Kneser graphs was
studied and the local chromatic numbers of the Kneser graphs
$K(2k+1,k)$ and ${\rm KG}(n,2)$ were determined. Also, it was
shown that for any positive integers $n$ and $k$ with $n\geq 2k$,
$\chi_l({\rm KG}(n,k))\leq 2\chi({\rm KG}(n,k))-2=2n-4k+2$. To
see this, for any $1\leq i\leq n-2k+1$ set
$$C_{2i-1}\isdef \{A\in V({\rm KG}(n,k)): \{1,2,\ldots ,i\}\cap A=\{i\}\ \},$$
$$C_{2n-4k+2}\isdef \{A\in V({\rm KG}(n,k)):  \{1,2,\ldots ,n-2k+1\}\cap A=\emptyset \}.$$
Now, it is readily seen that the aforementioned partition is a
local coloring, consequently, $\chi_s({\rm KG}(n,k))\leq
\chi_l({\rm KG}(n,k))\leq 2\chi({\rm KG}(n,k))-2=2n-4k+2$.

In \cite{OMPO1}, it was conjectured that for any positive
integers $n\geq 2k \geq 4$ we have $\chi_l({\rm
KG}(n,k))=2\chi({\rm KG}(n,k))-2=2n-4k+2$.

\begin{alphcon}\label{cc} {\rm \cite{OMPO1}}
For every positive integers $n$ and $k$ with $n\geq 2k\geq 4$, we
have $\chi_l({\rm KG}(n,k))=2\chi({\rm KG}(n,k))-2=2n-4k+2$.
\end{alphcon}

In view of Lemma \ref{BOUND}, we obtain $\chi({\rm KG}(n,k)) \geq
{\chi_l({\rm KG}(n,k))\over 2}+1$. Hence, the aforementioned
conjecture can be considered as a generalization of Kneser's
conjecture.

\begin{alphthm}{\rm \cite{OMPO1}}\label{MOO}
Let $n$ and $k$ be positive integers where $n\geq 2k \geq 4$. If
${\rm KG}(n,k)$ has a minimum local coloring with a color class
of size at least ${n-1 \choose k-1}-{n-k-1 \choose k-1}+2$, then
$\chi_l({\rm KG}(n,k))=\chi_l({\rm KG}(n-1,k))+2$.
\end{alphthm}

In the proof of the aforementioned theorem, the authors do not
employ the triangle condition (no triangle may be colored with
three consecutive numbers). Consequently, the above theorem holds
for star-free chromatic number as well. Here, we sketch the proof
of Theorem \ref{MOO}.  Let $c: V({\rm KG}(n, k))\rightarrow
\{1,2,\ldots ,\chi_s({\rm KG}(n,k))\}$ be a star-free coloring.
Set $C_i\isdef c^{-1}(i)$ for any $1\leq i \leq \chi_s({\rm
KG}(n,k))$. Also, assume that $C_j$ is a color class with size at
least ${n-1 \choose k-1}-{n-k-1 \choose k-1}+2$. The Hilton and
Milner theorem \cite{HIMI} says that if $I$ is an independent set
of ${\rm KG}(n, k)$ of size at least ${n-1 \choose k-1}-{n-k-1
\choose k-1}+2,$ then for some $i \in [n]$ we have $\cap_{A \in
I} A = \{i\}$. Hence, there exists some $i \in [n]$ such that
$\displaystyle \cap_{A \in C_j} A = \{i\}$. Now it is easy to
check that if $j\geq 2$ (resp. $j\leq n-1$), then $\displaystyle
i \in \cap_{A \in C_{j-1}} A $ (resp. $\displaystyle i\in \cap_{A
\in C_{j+1}} A $). By the above claim the vertices with the
colors $j-1 ,j ,j+1$  induce an empty subgraph. Without loss of
generality, suppose $i=n$. Now we define the coloring $c': V({\rm
KG}(n-1,k)) \rightarrow \{1,2,\ldots ,\chi_s({\rm KG}(n,k))-2 \}$
as follows.
$$c'(A)\isdef \left \{ \begin{array}{ll}
c(A) & if\ c(A) \leq j-2\\
c(A)-2 & if\ c(A) \geq j+2.
\end{array}\right. $$
One can check that $c'$ is a star-free coloring for ${\rm
KG}(n-1,k)$, hence, $\chi_s({\rm KG}(n-1,k))\leq \chi_s({\rm
KG}(n,k))-2$. On the other hand, any star-free coloring of ${\rm
KG}(n-1,k)$ with $t$ colors can be extended to a star-free
coloring of ${\rm KG}(n,k)$ with $t+2$ colors. Consequently,
$\chi_s({\rm KG}(n,k))=\chi_s({\rm KG}(n-1,k))+2$.

\begin{cor}\label{low}
Let $n\geq 2k \geq 4$ be positive integers. If ${\rm KG}(n,k)$
has a minimum star-free coloring with a color class of size at
least ${n-1 \choose k-1}-{n-k-1 \choose k-1}+2$, then $\chi_s({\rm
KG}(n,k))=\chi_s({\rm KG}(n-1,k))+2$.
\end{cor}

Section $2$ presents some preliminaries. In Section $3$, in view
of Tucker-Ky Fan's lemma, we show that for any Kneser graph ${\rm
KG}(n,k)$ we have $\chi_s({\rm KG}(n,k))\geq \max\{ 2\chi({\rm
KG}(n,k))-10, \chi({\rm KG}(n,k))\}$ where $n\geq 2k \geq 4$.
Moreover, we show that if $k$ is sufficiently large then
$\chi_l({\rm KG}(2k+t,k))=\chi_s({\rm KG}(2k+t,k))=2\chi({\rm
KG}(2k+t,k))-2=2t+2$. This gives a partial answer to Conjecture
\ref{cc}.
\section{Borsuk-Ulam Theorem and its Generalizations}
The Lov\'{a}sz's proof \cite{LOV} of the Kneser's conjecture was
an outstanding application of topological methods in 1978. In
fact, the proof of Lov\'{a}sz is based on the Borsuk-Ulam
theorem. The Borsuk-Ulam theorem says that if $f: S^n \rightarrow
\mathbb{R}^n$ is a continuous mapping from the unit sphere in
$\mathbb{R}^{n+1}$ into $\mathbb{R}^{n}$, then there exists a
point $x\in S^n$ where $f(x) = f(-x)$, that is, some pair of
antipodal points has the same image. There are several different
equivalent versions, various different proofs, several extensions
and generalizations, and many interesting applications for the
Borsuk-Ulam theorem that mark it as a great theorem, see
\cite{FAN, MAT1, MAT2, MEU, PRSU, SCH, SITA1, SITA2, TUC, SPE,
SPSU}.

Let $S^n$ denote the $n$-sphere, i.e.,
$S^n=\{x\in \mathbb{R}^{n+1}|\ ||x||= 1\}$. Assume that $T$ is a triangulation of $S^n$.
The triangulation $T$ is termed antipodally symmetric around the origin,
if $\sigma$ is a simplex in $T$, then $-\sigma$ is also
 a simplex in $T$. Tucker's lemma is a combinatorial analogue of the
Borsuk-Ulam theorem with several useful applications, see
\cite{MAT1, MAT2}.

\begin{alphlem}
{\rm(}Tucker's lemma. {\rm \cite{TUC} )} Let $T$ be a symmetric
triangulation of the $n$-sphere $S^n$ where $n$ is a positive
integer. Assume that each vertex $u$ of $T$ is assigned a label
$\lambda(u)\in \{\pm 1,\pm 2,\ldots ,\pm n\}$ such that $\lambda$
is an antipodal map, i.e., $\lambda(-u)=-\lambda(u)$ for any
vertex $u$ of $T$. Then some pair of adjacent vertices of $T$
have labels that sum to zero.
\end{alphlem}

Another interesting generalization of Borsuk-Ulam theorem is Ky
Fan's lemma \cite{FAN}, which generalizes the
Lusternik–Schnirelmann theorem which is the version of the
Borsuk-Ulam theorem involving a cover of the $n$-sphere by $n+1$
sets, all open or all closed. Just like the Borsuk-Ulam theorem
it has several equivalent forms, see \cite{FAN}.

\begin{alphlem}
{\rm (}Ky Fan's lemma.{\rm )} Let $n$ and $k$ be two arbitrary
positive integers. Assume that $k$ closed subsets {\rm (}resp.
open subsets{\rm )} $F_1, F_2,\ldots, F_k$ of the $n$-sphere
$S^n$ cover $S^n$ and also no one of them contains a pair of
antipodal points. Then there exist $n+2$ indices $l_1, l_2,
\ldots,l_{n+2}$, such that $1\leq l_1 < l_2 < \cdots <
l_{n+2}\leq k$ and
$$F_{l_1}\cap -F_{l_2}\cap
F_{l_3}\cap \cdots \cap (-1)^{n+1}F_{l_{n+2}}\not = \emptyset,.$$
where $-F_i$ denotes the antipodal set of $F_i$. In particular,
$k \geq n+2$.
\end{alphlem}

This lemma has useful applications in graph colorings and
provides useful information about coloring properties of Kneser
graphs, see  \cite{MEU, SITA1, SITA2, ZHU}. For instance, it was
shown that the circular chromatic number and the chromatic number
of the Kneser graph ${\rm KG}(n,k)$ are equal provided that $n$
is even, see \cite{MEU,SITA1}. Moreover, we can consider the
subcoloring theorem as an interesting application of Ky Fan's
lemma, see \cite{SPE, SPSU}. This theorem states if $c$ is a
proper coloring of the Kneser graph ${\rm KG}(n,k)$ with $m$
colors, then there exists a multicolored complete bipartite graph
$K_{{\lceil{ r\over 2}\rceil},{\lfloor{ r\over 2}\rfloor}}$ with
$r\isdef \chi({\rm KG}(n,k))$ such that the $r$ different colors
occur alternating on the two sides of the bipartite graph with
respect to their natural order. The subcoloring theorem has been
generalized for general Kneser graphs in \cite{SITA2}.

In \cite{FAN}, Fan introduced a generalization of Tucker's lemma
which is called  Tucker-Ky Fan's lemma.

\begin{alphlem}{\rm (}Tucker-Ky Fan's lemma.{\rm )}
Let $T$ be a symmetric triangulation of $S^n$ and $m$ be a fixed
positive integer. Also, assume that each vertex $u$ of $T$ is
assigned a label $\lambda(u)\in \{\pm 1,\pm 2,\ldots ,\pm m\}$
such that $\lambda$ is an antipodal map and furthermore labels at
adjacent vertices do not sum to zero. Then there are an odd
number of $n$- simplices whose labels are of the form $\{a_0,-
a_1, ..., (- 1)^na_n\}$, where  $1\leq a_0 < a_1 < \cdots < a_n
\leq m$. In particular $m \geq n + 1$.
\end{alphlem}

Now, we introduce a special triangulation of $S^{n-1}$.  Let $e_1,
e_2,\ldots , e_n$ be the vectors of the standard orthonormal
basis of $R^n$ ($e_i$ has a $1$ at position $i$ and $0$'s
elsewhere). Define a simplicial complex $C^{n-1}$  (cross
polytope) as follows. The vertex set of $C^{n-1}$ is equal to
$\{\pm e_1,\pm e_2,\ldots ,\pm e_n\}$. Also, a subset $F\subseteq
\{\pm e_1,\pm e_2,\ldots ,\pm e_n\}$  forms the vertex set of a
face of the cross polytope if there is no $i \in [n]$ with both
$e_i\in F$ and $-e_i\in F$.

Set $V_n \isdef \{- 1,0,1\}^n$. Consider a partial ordering $\leq$
on $V_n$ that for any $u,v \in V_n$ we have $u \leq v$ if $u_i
\leq v_i$ for any $1\leq i \leq n$ where $0 \leq -1$ and $0 \leq
1$. Now, we introduce a symmetric triangulation of $S^{n-1}$, say $K$, whose vertex set can be identified with the
vectors of $V_n\setminus \{(0, 0,\ldots,0)\}$.

Let $K$ be the the first barycentric subdivision of $C^{n-1}$. Thus,
the vertices of $K$ are centers of gravity of the simplices of
$C^{n-1}$ and the simplices of $K$ correspond to chains of
simplices of $C^{n-1}$ under inclusion. In fact, the vertex set of $K$ can be identified with
$V_n\setminus \{(0, 0,\ldots,0)\}$. Moreover, a simplex of $K$ can be
identified uniquely with a chain in the set $V_n\setminus \{(
0,0,\ldots,0)\}$ under $\leq$. Also, one can see that $K$ is a symmetric
triangulation of $S^{n-1}$ which meets conditions in Tucker-Ky Fan's lemma. Hence, to
use the Tucker-Ky Fan's lemma, it is sufficient to define a
labeling $\lambda: V_n\setminus \{(0, 0,\ldots,0)\} \longrightarrow \{\pm 1, \pm 2, \ldots,
\pm m\}$ which has the antipodal property, i.e., for any nonzero
$v\in V_n$ we have $\lambda(-v)=-\lambda(v)$, and furthermore,
for any $u,v \in V_n\setminus \{(0, 0,\ldots,0)\}$ where $u\leq v$ or $v\leq u$, we have
$\lambda(u)+\lambda(v)\not =0$, that is, labels at adjacent
vertices do not sum to zero.
\section{Star-Free Chromatic Number of Kneser Graphs}
It was shown in \cite{OMPO1} if $n\geq 2k^3-2k^2$, then
$\chi_l({\rm KG}(n,k))=\chi_l({\rm KG}(n-1,k))+2$.

\begin{thm}\label{k2}
Let $n\geq 2k \geq 4$ be positive integers where $n\geq
2k^3-2k^2-2k+4$. Then $\chi_s({\rm KG}(n,k))=\chi_s({\rm
KG}(n-1,k))+2$.
\end{thm}
\begin{proof}{Let $c$ be a minimum star-free coloring of ${\rm
KG}(n,k)$. In view of Corollary \ref{low}, it is sufficient to
show that there exists a color class of $c$ with size at least
${n-1 \choose k-1}-{n-k-1 \choose k-1}+2$. Note that $\chi_s({\rm
KG}(n,k))\leq 2 \chi({\rm KG}(n,k))-2$ and also it was shown in
\cite{HAZH1} that if $n\geq 2k^2(k-1)$, then ${{n \choose k}
\over 2(n-2k+2)}\geq {n-1 \choose k-1}-{n-k-1 \choose k-1}+2$.
Hence, for $k=2$ the assertion follows. Now assume that $k\geq
3$. By double counting we have
\begin{equation}\label{eq1}
{n-1 \choose k-1}-{n-k-1 \choose k-1}+2\leq k{n-2 \choose k-2}.
\end{equation}
Hence, if we show that there exists a color class of $c$ with
size at least $k{n-2 \choose k-2}$, then $\chi_s({\rm
KG}(n,k))=\chi_s({\rm KG}(n-1,k))+2$. For any positive integers
$a$ and $b$ where $a\geq 2k-2$ and $b\leq 2k-3$ we have
\begin{equation}\label{eq2}
a(a-b-2k+1) \leq (a-{b+2k-1 \over 2})^2 \leq (a-b)(a-b-1).
\end{equation}
Let $b\leq 2k-3$ and $a\geq 2k^3-2k^2$. Set $n\isdef a-b$. In
view of (\ref{eq2}), it is readily seen that if $n\geq
2k^3-2k^2-2k+3$, then
$$(2k^3-2k^2)(n-2k+1)\leq a(a-b-2k+1) \leq (a-b)(a-b-1)=n(n-1).$$
Consequently,
$${{n \choose k} \over k{n-2 \choose k-2}}={n(n-1) \over k^2(k-1)}
\geq 2(n-2k+1)=2\chi({\rm KG}(n,k))-2 \geq \chi_s({\rm
KG}(n,k)).$$ Therefore, there is a color class of $c$ with size
at least $k{n-2 \choose k-2}$, as desired.}
\end{proof}

\begin{rem}
In \cite{ALHA}, the inequality {\rm (\ref{eq1})} has been applied
to show that if $n\geq 2k^3-2k^2-2k+4$, then the circular
chromatic number and the chromatic number of the Kneser graph
${\rm KG}(n,k)$ are equal  {\rm (}for definition and more on
circular chromatic number see {\rm \cite{ZHU})}.
\end{rem}

Matou{\v{s}}ek \cite{MAT2} introduced an interesting proof for
the Kneser conjecture by using Tucker's lemma. Similarly, in view
of Tucker-Ky Fan's lemma, one can show that for any positive
integers $n\geq 2k$, we have $\chi_s({\rm KG}(n,k))\geq 2\chi({\rm
KG}(n,k))-10=2n-4k-6$.
\begin{thm}\label{main}
For every positive integers $n\geq 2k\geq 4$, we have
$\chi_s({\rm KG}(n,k))\geq \max \{ 2\chi({\rm KG}(n,k))-10,
\chi({\rm KG}(n,k))\}.$
\end{thm}
\begin{proof}{
Let $k$ be a fixed positive integer. In view of Theorem \ref{k2},
it is sufficient to show that $\chi_s({\rm KG}(n,k))\geq
2\chi({\rm KG}(n,k))-10$ provided that $n< 2k^3-2k^2-2k+4$. Thus,
assume that $2k\leq n\leq  2k^3-2k^2-2k+3$. Also, on the
contrary, suppose that $\chi_s({\rm KG}(n,k))\leq 2\chi({\rm
KG}(n,k))-11=2n-4k-7.$

Suppose that $c$ is a star-free coloring of the Kneser graph
${\rm KG}(n,k)$ with $2n - 4k -7$ colors. Suppose that the colors
are numbered $2k-1, 2k,\ldots , m=2n-2k-9$. Now, we introduce a
labeling $\lambda: V_n\setminus \{(0, 0,\ldots,0)\} \longrightarrow \{\pm 1, \pm 2, \ldots,
\pm m\}$. Consider an arbitrary linear ordering $\leq$ on power
set of $[n]$ that refines the partial ordering according to size,
that is, if $|A|<|B|$ then $A<B$.

Let $w=(w_1,w_2,\ldots,w_n) \in V_n\setminus \{(0, 0,\ldots,0)\}$. To define $\lambda(w)$, we
consider the ordered pair $(P(w),N(w))$ of disjoint subsets of
$[n]$ defined by
$$P(w)\isdef \{i\in [n]: w_i=+1\}\ {\rm and}\ N(w)\isdef \{i\in [n]: w_i=-1\}.$$
We consider two cases. If $|P(w)|+|N(w)|\leq 2k-2$ (Case I) then
we set
$$\lambda(w)\isdef \left \{ \begin{array}{ll}
|P(w)|+|N(w)| & if\ P(w)\geq N(w)\\
-|P(w)|-|N(w)| & if\ P(w)< N(w).
\end{array}\right. $$
Now assume that  $|P(w)|+|N(w)|\geq 2k-1$ (Case II). Note that if
$|P(w)|+|N(w)|\geq 2k-1$ then at least one of $P(w)$ and $N(w)$
has size at least $k$. If $P(w) \geq N(w)$ (resp. $P(w) < N(w)$)
we define $\lambda(w)=t$ (resp. $\lambda(w)=-t$) where $t$ is the
largest positive integer such that there exist two distinct
$k$-subsets $A, B\subseteq P(w)$ (resp. $A, B\subseteq N(w)$)
where $c(A)=c(B)=t$, otherwise, $t$ is the largest positive
integer such that there exists a $k$-subset $A\subseteq P(w)$
(resp. $A\subseteq N(w)$) where $c(A)=t$.

It is a simple matter to check that the labeling $\lambda$ is
well-defined and that it has the antipodal property, i.e., for
any nonzero $v\in V_n$ we have $\lambda(-v)=-\lambda(v)$.
Furthermore, labels at adjacent vertices do not sum to zero.
Hence, in view of Tucker-Ky Fan's lemma, there exist an $(n-1)$-
simplex (a chain of length $n$ in $V_n$), say $\sigma$, whose
labels are of the form $\{a_0,- a_1, ..., (- 1)^{n-1}a_{n-1}\}$,
where $1\leq a_0 < a_1 < \cdots < a_{n-1} \leq m= 2n-2k-9$. Let
$V(\sigma)=\{v_1,v_2,\ldots,v_{n}\}$. Referring to our
construction of $\lambda$, at least the
$n-(2k-2)=n-2k+2=\chi({\rm KG}(n,k))$ highest of these labels
were assigned by Case II. In view of Tucker-Ky Fan's lemma, as
$m=2n-2k-9$, we can deduce that there exist at least $2k+8$ pairs
of vertices of $\sigma$ such that for any pair $\{v_i, v_j\}$ we
have $|\lambda(v_i)+\lambda(v_j)|=1$. Therefore, there exist two
vertices $v_i, v_j \in V(\sigma)$ with $|\lambda(v_i)|\geq 2k-1$,
$|\lambda(v_j)|\geq 2k-1 $, and $|\lambda(v_i)+\lambda(v_j)|=1$
such that $|P(v_i)|+|N(v_j)| \geq 2k+9$ or $|P(v_j)|+|N(v_i)|
\geq 2k+9$. Without loss of generality, assume that $P(v_i)
> N(v_i)$, and consequently, $|P(v_i)|+|N(v_j)| \geq 2k+9$. In
view of definition of $\lambda$ and since $c$ is a star-free
coloring and $|\lambda(v_i)+\lambda(v_j)|=1$, one can conclude
that all of $k$-subsets of $P(v_i)$ and $N(v_j)$ receive distinct
colors. Hence, considering $n\leq 2k^3-2k^2-2k+3$, we should have
$${|P(v_i)|\choose k}+{|N(v_j)|\choose k} \leq 2\chi({\rm KG}(n,k))-11\leq 4k^3-4k^2-8k-1.$$
On the other hand, $|P(v_i)|+|N(v_j)| \geq 2k+9$ so
$${k+5 \choose k}+{k+4 \choose k}\leq  {|P(v_i)|\choose k}+{|N(v_j)|\choose
k}.$$ But, one can check that for any $k\geq 2$, ${k+5 \choose
k}+{k+4 \choose k} > 4k^3-4k^2-8k-1$ which is a contradiction.

}
\end{proof}
\begin{rem}
Similarly, one can show that if $n\geq 2k$ and $k \geq 81$, then
$\chi_s({\rm KG}(n,k))\geq \max \{ 2\chi({\rm KG}(n,k))-8,
\chi({\rm KG}(n,k))\}$.
\end{rem}

It seems that the star-free chromatic number of the Kneser graph
${\rm KG}(n,k)$ is equal to $2\chi({\rm KG}(n,k))-2$ provided
that $n\geq 2k\geq 4$.

\begin{thm}\label{3k}
For any positive integers $n \geq 2k \geq 4$, if $n \leq {8\over
3}k$, then $\chi_s({\rm KG}(n,k))=2\chi({\rm KG}(n,k))-2=2n-4k+2$.
\end{thm}
\begin{proof}{On the contrary,
let $\chi_s({\rm KG}(n,k))\leq 2\chi({\rm KG}(n,k))-3=2n-4k+1$
provided that $n\leq {8\over 3}k$. The proof is almost similar to
that of Theorem \ref{main}. The labeling $\lambda$ and
$(n-1)$-simplex $\sigma$ are defined similarly. As $m= 2n-2k-1$,
we can deduce that there exist at least $2k$ pairs of vertices of
$\sigma$ such that for any pair $\{v_i, v_j\}$ we have
$|\lambda(v_i)+\lambda(v_j)|=1$. Therefore, there exist two
vertices $v_i$ and $v_j$ with $|\lambda(v_i)|\geq 2k-1$,
$|\lambda(v_j)|\geq 2k-1$, and $|\lambda(v_i)+\lambda(v_j)|=1$
such that $|P(v_i)|+|N(v_j)| \geq 2k+1$ or $|P(v_j)|+|N(v_i)|
\geq 2k+1$. Without loss of generality, assume that $P(v_i)
> N(v_i)$, and consequently, $|P(v_i)|+|N(v_j)| \geq 2k+1$. In
view of definition of $\lambda$ and since $c$ is a star-free
coloring, one can conclude that all of $k$-subsets of $P(v_i)$
and $N(v_j)$ receive distinct colors. Without loss of generality,
suppose that $|P(v_i)|\geq k+1$. For any vertex $v_j\in \sigma$ we
have $P(v_i)\cap N(v_j)=\emptyset$. Consequently, in view of
Tucker-Ky Fan's lemma, at least $\lfloor {\chi({\rm
KG}(n,k))\over 2}\rfloor$ colors do not assign to the vertices of
${P(v_i)\choose k}$ which implies that
$${|P(v_i)|\choose k}+ {\chi({\rm KG}(n,k))-1\over 2}\leq 2\chi({\rm KG}(n,k))-3.$$
Therefore,
$$k+1\leq {3 \over 2}\chi({\rm KG}(n,k))-{5\over 2}={3 \over 2}(n-2k+2)-{5 \over 2},$$
consequently, $n \geq {8\over 3}k+{1\over 3}$ which is a contradiction.}
\end{proof}

We know that $\chi_s(G)\leq \chi_l(G)$. Hence, we have the
following corollary.

\begin{cor}
For any positive integers $n \geq 2k \geq 4$, if $n \leq {8\over
3}k$, then $\chi_l({\rm KG}(n,k))=2\chi({\rm KG}(n,k))-2=2n-4k+2$.
\end{cor}

Theorems \ref{main} and \ref{3k} motivate us to propose the
following conjecture which can be considered as a generalization
of Kneser's conjecture.

\begin{con}\label{HHH}
For any positive integers $n\geq 2k\geq 4$, we have $\chi_s({\rm
KG}(n,k))=2\chi({\rm KG}(n,k))-2=2n-4k+2$.
\end{con}
\bibliographystyle{plain}
\bibliography{xbbl}

\end{document}